\documentclass[11pt]{article}
\usepackage{amsthm,amssymb,amscd,amsmath,latexsym,indentfirst,color,amsfonts}
\usepackage[mathcal]{eucal}
\theoremstyle{plain}
\newtheorem{theorem}{Theorem}[section]
\newtheorem{prop}[theorem]{Proposition}

\newtheorem{lem}[theorem]{Lemma}
\theoremstyle{remark}
\newtheorem{rem}[theorem]{Remark}

\def\g{\frak g}\def\h{\frak h}
\def\a{\frak a}
   
\def\b{\frak b}\def\n{\frak n}
 \def\ld{\lambda}
\def\al{\alpha}\def\ga{\gamma}
\def\w{\wedge}\def\rt{\rightarrow}\def\st{\subset}\def\ddt{\Delta}

\def\Ga{\Gamma}\def\ga{\gamma}
\newcommand{\be}{\begin {equation}}
\newcommand{\ee}{\end {equation}}
\newcommand{\bp}{\begin {proof}}
\newcommand{\ep}{\end {proof}}
\newcommand{\bee}{\begin {equation*}}
\newcommand{\eee}{\end {equation*}}
\newcommand{\lb}{\label}

\begin{document}

\title{Maximal eigenvalues of a Casimir operator and multiplicity-free modules
\footnote{Research supported by NSFC Grant No.10801116 and by 'the
Fundamental Research Funds for the Central Universities'} }
\author{Gang HAN \\
Department of Mathematics, College of Science,\\Zhejiang University,\\{Hangzhou 310027, China}\\
{E-mail: mathhg@hotmail.com}\\Telephone: 086-0571-87953843\\Fax:
086-0571-87953794 }
\date{January 15, 2011}

\maketitle

\begin{abstract}
Let $\g$ be a finite-dimensional complex semisimple Lie algebra and
$\b$ a Borel subalgebra. Then $\g$ acts on its exterior algebra
$\w\g$ naturally.
 We prove that the maximal
eigenvalue of the Casimir operator on $\w\g$ is one third of the
dimension of $\g$, that the maximal eigenvalue $m_i$ of the Casimir
operator on $\w^i\g$ is increasing for $0\le i\le r$, where $r$ is
the number of positive roots, and that the corresponding eigenspace
$M_i$ is a multiplicity-free $\g$-module whose highest weight
vectors corresponding to certain ad-nilpotent ideals of $\b$. We
also obtain a result describing the set of weights of the
irreducible representation of $\g$ with highest weight a multiple of
$\rho$, where $\rho$ is one half the sum of positive roots.

\end{abstract}

%\textbf{2010 Mathematics Subject Classification}:17B10\bigskip

%\textbf{Key words}: Casimir operator, exterior algebra, multiplicity-free module

 \section{Introduction}
 \setcounter{equation}{0}\setcounter{theorem}{0}
Let $\g$ be a finite-dimensional complex semisimple Lie algebra and
$U(\g)$ its universal enveloping algebra. The study of the
$\g$-module structure of its exterior algebra $\w\g$ has a long
history. Although this module structure is still not fully
understood, Kostant
 have done a lot of important work on it, see for example \cite{k1} and \cite{k2}.

Let $Cas\in U(\g)$ be the Casimir element with respect to the
Killing form. Let $m_i$ be the maximal eigenvalue of $Cas$ on
$\wedge ^{i} \g$ and $M_i$ be the corresponding eigenspace. Let $p$
be the maximal dimension of commutative subalgebras of $\g$. In
\cite{k1} it is proved that $m_i\le i$ for any $i$ and $m_i=i$ for
$0\le i\le p$, and if $m_i=i$ then $M_i$ is a multiplicity-free
$\g$-module whose highest weight vectors corresponding to
$i$-dimensional abelian ideals of $\b$. The integer $p$ for all the
simple Lie algebras was determined by Malcev, and Suter
 gave a uniform formula for $p$ in \cite{s}.

Fix a Cartan subalgebra $\h$ of $\g$ and a set $\ddt^+$ of positive
roots. Let $\rho\in \h^{*}$ be one half the sum of all the positive
roots. For any $\ld\in\h^{*}$, let $V_{\ld}$ denote the irreducible
representation of $\g$ with highest weight $\ld$.

In this paper we will prove the following result, which extends some theorems of Kostant. Let $n=dim \ \g$ and $r$ be the number of positive
roots.
\begin{theorem} [Theorem \ref{e}]

One has $m_i\le n/3$ for $i=0,1,\cdots,n$, and $m_i=n/3$ if and only if
$i=r,r+1,\cdots,r+l$. For $s=0,1,\cdots,l$, $M_{r+s}=\left(
             \begin{array}{c}
               l \\
               s \\
             \end{array}
           \right) V_{2\rho}$.

For $0\le i< r$ one has $m_i< m_{i+1}$. For $1\le i\le r$, $M_i$ is a multiplicity-free $\g$-module, whose highest
weight vectors corresponding to certain ad-nilpotent ideals of $\b$. In fact $\oplus_{i=0}^r M_i$ is also a
multiplicity-free $\g$-module.
\end{theorem}
This result relates $M_i$ to ad-nilpotent ideals of $\b$, which are
classified in \cite{pp}. But it will be complicated to
determine those ad-nilpotent ideals of $\b$ corresponding to the
highest weight vectors of $M_i$.

To prove this theorem, we need the following interesting result.

\begin{prop}[Proposition \ref{d}]
Let $k\in \mathbb{Z}^+$. The set of weights of $V_{k\rho}$ (whose
dimension is $(k+1)^r$) is $$\{\sum_{i=1}^r c_i\al_i|\al_i\in\ddt^+,
c_i=-k/2,-k/2+1,\cdots,k/2-1,k/2\}.$$
\end{prop}

\section{Weights of a representation with highest weight a multiple of $\rho$}
 \setcounter{equation}{0}\setcounter{theorem}{0}
Let $\g$ be a finite-dimensional complex semisimple Lie algebra. Fix
a Cartan subalgebra $\h$ of $\g$ and a Borel subalgebra $\b$ of $\g$
containing $\h$. Let $\ddt$ be the set of roots of $\g$ with respect
to $\h$ and $\ddt^+$ be the set of positive roots whose
corresponding root spaces lie in $\b$. Let $W$ be the Weyl group.
Let $\n=[\b,\b]$. Then $\b=\h\oplus\n$. An ideal of $\b$ contained
in $\n$ is called an ad-nilpotent ideal of $\b$, as it consists of
ad-nilpotent elements. Let $\Gamma\st \h^{*}$ be the lattice of
$\g$-integral linear forms on $\h$ and $\Lambda\st\Gamma$ the subset
of dominant integral linear forms. Let $(,)$ be the bilinear form on
$\h^{*}$ induced by the Killing form. Let $l=dim\ \h$, $r=|\ddt^+|$
and $n=l+2r=dim\ \g$. Assume $\ddt^+=\{\al_1,\al_2,\cdots,\al_r\}.$

 For any $\ld\in \Lambda$,
let $\pi_\ld:\g\rt End(V_\ld)$ be the irreducible representation of
$\g$ with highest weight $\ld$, and $\Ga(V_\ld)$ be the set of
weights, with multiplicities. Any $\ga\in\Ga(V_\ld)$ will appear $k$
times if the dimension of the $\ga$-weight space is $k$. For example
$\Ga(\g)=\ddt\cup\{0,\cdots,0\}$ ($l$ times). If $U\st V_\ld$ is an
$\h$-invariant
 subspace then we will also use $\Ga(U)$ to denote the the set of weights of $U$ with multiplicities and define $$<U>=\sum_{\ga\in\Ga(U)} \ga.$$
 For any $S\st \Ga(V_\ld)$, we also define $<S>=\sum_{\ga\in S} \ga.$

  Let $\rho\in
\h^{*}$ be one half the sum of all the positive roots. For any $k\in
\mathbb{Z}^+$, the representation $V_{k\rho}$ of $\g$ has dimension
$(k+1)^r$ by Weyl's dimension formula. The following result describes the set of weights of
$V_{k\rho}$, which is well-known if $k=1$ (see e.g. \cite{k0}).
\begin{prop}\lb{d}
The set of weights of $V_{k\rho}$ is $$\Ga(V_{k\rho})=\{\sum_{i=1}^r
c_i\al_i|\al_i\in\ddt^+, c_i=-k/2,-k/2+1,\cdots,k/2-1,k/2\},$$ or
equivalently,
$$\Gamma(V_{k\rho})=\{k\rho-\sum_{i=1}^r c_i\al_i|\al_i\in\ddt^+,
c_i=0,1,\cdots,k.\}.$$
\end{prop}
\bp By Weyl's denominator formula
$$\prod_{i=1}^r (e^{\frac{k+1}{2}\al_i}-e^{-\frac{k+1}{2}\al_i})=\sum_{w\in W} sgn(w)
e^{w((k+1)\rho)}.$$ Then for $c_i=-k/2,-k/2+1,\cdots,k/2-1,k/2$ with
$i=1,\cdots,r$,

\bee
\begin{split}
\sum_{c_1,\cdots,c_r} e^{\sum_{i=1}^r
c_i\al_i}&=\prod_{i=1}^r(e^{(-\frac{k}{2})\al_i}+e^{(-\frac{k}{2}+1)\al_i}+\cdots+e^{(\frac{k}{2}-1)\al_i}+e^{(\frac{k}{2})\al_i})\\
&=\prod_{i=1}^r\frac{e^{\frac{k+1}{2}\al_i}-e^{-\frac{k+1}{2}\al_i}}{e^{\frac{1}{2}\al_i}-e^{-\frac{1}{2}\al_i}}\\&=\frac{\sum_{w\in
W} \ sgn(w) e^{w((k+1)\rho)}}{\prod_{i=1}^r
(e^{\frac{1}{2}\al_i}-e^{-\frac{1}{2}\al_i})}\\&=char(V_{k\rho}).
\end{split}
\eee  \ep

 Let $Cas\in U(\g)$ be the Casimir element corresponding
to the Killing form. For any $\ld\in\Ga$, define
$$Cas(\ld)=(\ld+\rho,\ld+\rho)-(\rho,\rho).$$ The following result is well-known.
\begin{lem}\lb{c}
If $\ld\in\Lambda$ then $Cas(\ld)$ is the scalar value
taken by $Cas$ on $V_\ld$. For any $\mu\in \Ga(\ld)$ one has $Cas(\mu)\le Cas(\ld)$ and $Cas(\mu)<Cas(\ld)$ if $\mu\neq\ld$.
\end{lem}
\section{Maximal eigenvalues of a Casimir operator and the corresponding eigenspaces}
 \setcounter{equation}{0}\setcounter{theorem}{0}

Let $\wedge \g$ be the exterior algebra of $\g$. Then $\g$ acts on
$\wedge \g$ naturally. Let $m_i$ be the maximal eigenvalue of $Cas$
on $\wedge ^{i} \g$ and $M_i$ be the corresponding eigenspace.

 One knows that $\wedge ^{i}\g$ is isomorphic to $\w
^{n-i}\g$ as $\g$-modules for each $i$, so one has
$$m_i=m_{n-i}$$ and $$M_i\cong M_{n-i}.$$ Let $p$ be the maximal dimension of commutative
subalgebras of $\g$. Kostant showed that $m_i\le i$ and $m_i=i$ for
$0\le i\le p$, and if $m_i=i$ then $M_i$ is spanned by $\wedge ^{k}
\a$, where $\a$ runs through $k$-dimensional commutative subalgebras
of $\g$.

 A nonzero vector $w\in\w\g$ is called decomposable if
$w=z_1\w z_2\w\cdots\w z_k$ for some
 positive integer $k$, where $z_i\in \g$. In this case let
$\a(w)$ be the respective $k$-dimensional subspace spanned by
$z_1,z_2,\cdots,z_k$.
\begin{theorem} [Proposition 6 and Theorem 7 of \cite{k1}]\lb{b}
(1) Let $$w=z_1\w z_2\w\cdots\w z_k\in\w^k\g$$ be a decomposable
vector. Then $w$ is a highest weight vector if and only if $\a(w)$
is $\b$-normal, i.e., $[\b,\a(w)]\st\a(w)$. In this case the highest
weight of the simple $\g$-module generated by $w$ is $<\a(w)>$.

 Thus there is a one-to-one correspondence between all the decomposably-generated simple $\g$-submodules of $\w^k\g$ and all the
  $k$-dimensional $\b$-normal subspaces of $\g$.

(2) Let $\a_1,\a_2$ be any two ideals of $\b$ lying in $\n$. Then
$<\a_1>=<\a_2>$ if and only if $\a_1=\a_2$. Thus, if
$V_1\st\w^k\g,V_2\st\w^j\g$ are two decomposably-generated simple
$\g$-submodules which corresponds to ideals of $\b$ lying in $\n$,
then $V_1$ is equivalent to $V_2$ if and only if $V_1=V_2$.

\end{theorem}
\begin{theorem} \lb{e} (1) One has $$m_i=max\{||\rho+\ga_1+\cdots+\ga_i||^2-||\rho||^2\ |\{\ga_t|t=1,\cdots,i\}\st\Ga(\g)\}$$ for any $i$.

(2)One has $m_i\le n/3$ for $i=0,1,\cdots,n$, and $m_i=n/3$ if and
only if $i=r,r+1,\cdots,r+l$. For $s=0,1,\cdots,l$, $M_{r+s}=\left(
             \begin{array}{c}
               l \\
               s \\
             \end{array}
           \right) V_{2\rho}$.

(3)For $0\le k< r$ one has $m_k< m_{k+1}$. For $1\le k\le r$, $M_k$ is a multiplicity-free $\g$-module, whose highest
weight vectors corresponding to those $k$-dimensional ad-nilpotent ideals $\a$ of $\b$ such that $Cas(<\a>)=m_k$. In fact $\oplus_{k=0}^r M_k$ is also a
multiplicity-free $\g$-module.
\end{theorem}

\bp (1)  For $j=1,\cdots,r,$ let $x_j$ (resp. $y_j$) be a weight
vector corresponding to $\al_j$ (resp. $-\al_j$). Let
$\{h_1,\cdots,h_l\}$ be a basis of $\h$. Then
$$A=\{x_1,\cdots,x_r,y_1,\cdots,y_r,h_1,\cdots,h_l\}$$ is a basis of
$\g$ consisting of weight vectors. Then $$B_i=\{a_1\wedge a_2\wedge
\cdots\wedge a_i| a_j\in A\}$$ is a basis of $\w^{i}\g$ consisting
of weight vectors. Let $$C_i=\{v\in B_i\ | Cas(<\a(v)>)=m_i\}.$$
Then by Corollary 2.1 of \cite{k1} $M_i$ is the direct sum of simple
$\g$-modules with highest weight vectors $v\in C_i$. It is clear
that $$Cas(<\a(v)>)=||\rho+\ga_1+\cdots+\ga_i||^2-||\rho||^2$$ if
the weight of $a_j$ is $\ga_j$, thus (1) follows.

(2)For any $S=\{\ga_j|j=1,\cdots,i\}\st\Ga(\g)$, $<S>$ is a weight
of $\pi_{2\rho}$ by Proposition \ref{d}. Thus by Lemma \ref{c}
$Cas(<S>)\le Cas(2\rho)=8||\rho||^2=n/3$, as $||\rho||^2=n/24$. So
$m_i=n/3$ if and only if there exists $S\st\Ga(\g)$ such that
$|S|=i$ and $<S>=2\rho$. Then $S$ must be of the form
$\{x_1,\cdots,x_r,h_{j_1},\cdots,h_{j_s}\}$ and thus $r\le i\le
r+l$. For $0\le s\le l$, it is clear $$C_{r+s}=\{x_1\w\cdots\w x_r\w
h_{j_1}\w\cdots\w h_{j_s}|1\le j_1<j_2<\cdots<j_s\le l \},$$ thus
$M_{r+s}=\left(
             \begin{array}{c}
               l \\
               s \\
             \end{array}
           \right) V_{2\rho}$.

(3)We first show $m_{k+1}>m_{k}$ for $0\le k<r$, which clearly holds
in the case $k=0$. Assume $1\le k<r$. Let $v=a_1\w\cdots\w a_k\in
C_k$. Then $v$ is a highest weight vector of $M_k$, whose weight is
$<S>$ with $S=\Ga(\a(v))$. Then $Cas(<S>)=m_k$, and
$[\b,\a(v)]\st\a(v)$ by Theorem \ref{b} (1). Recall that for
$\ga=\sum_{i=1}^l k_i \ga_i\in\ddt^+$ where $\{\ga_i|i=1,\cdots,l\}$
is the set of simple roots, its height is defined as $\sum_{i=1}^l
k_i$. Choose a positive root $\al$ in $\ddt^+\setminus S$ (which is
nonempty as $k<r$) with largest height. Set $T=S\cup \{\al\}$. Let
$a\in A$ be the $\al$-weight vector and let $u=v\w a\in B_{k+1}$. By
the choice of $\al$ it is clear that $[\b,\a(u)]\st\a(u)$, thus $u$
is also a highest weight vector, whose weight is $<T>=<S>+\al$. As
$$(<T>,\al)=(<S>,\al)+(\al,\al)>0,$$ $<S>\in \Ga(V_\ld)$ with $\ld={<T>}$. Then
$$m_{k+1}\geq Cas(<T>)>Cas(<S>)=m_k.$$

Now assume $1\le k\le r$. Let $v=a_1\w\cdots\w a_k\in C_k$ , and let
$S=\Ga(\a(v))$. We will show $S\st\ddt^+$. If not, let
$S^{'}=S\setminus(S\cap(-S))$. Then $<S^{'}>=<S>$ and $|S^{'}|=t<k$.
Thus $m_k=Cas(S)=Cas(S^{'})\le m_t$, which contradicts to the
previous result. Thus for $1\le k\le r$ one always has $S\st\ddt^+$.

Any $v\in C_k$ is a highest weight vector, so $[\b,\a(v)]\st \a(v)$.
And if $1\le k\le r$ we have just showed $\Ga(\a(v))\st \ddt^+$.
Thus $\a(v)$ is an ad-nilpotent ideal of $\b$. Let $\ld(v)=<\a(v)>$.
Then
$$M_k=\oplus_{v\in C_k}V_{\ld(v)}.$$ By Theorem \ref{b} (2), if
$v_1,v_2\in C_k$ with $v_1\neq v_2$, then $\a(v_1)\neq\a(v_2)$ and
$\ld(v_1)\neq \ld(v_2)$. Thus $M_k$ is a multiplicity-free
$\g$-module, whose highest weight vectors corresponding to the
ad-nilpotent ideals $\a$ of $\b$ such that $Cas(<\a>)=m_k$. By
Theorem \ref{b} (2) one can further get that $\oplus_{k=0}^r M_k$ is
also a multiplicity-free $\g$-module.

\ep
\begin{rem}
Considering the isomorphism of $\g$-modules $\w^k\g$ and
$\w^{n-k}\g$, $\w^k\g$ is multiplicity-free for $0\le k\le r$ and
$n-r\le k\le n$. For $r\le k\le r+l$ ($r+l=n-r$), we have showed
that $M_k$ is primary of type $\pi_{2\rho}$. As $\g$-modules one has
$\w\g=2^l V_\rho\otimes V_\rho$ (see \cite{k2}), so $\w\g$ contains
exactly $2^l$ copies of $V_{2\rho}$, which is just $\oplus_{s=0}^l
M_{r+s}$.
\end{rem}

\end{document}